\documentclass[conference]{IEEEtran}
\usepackage{amsmath, mathrsfs,amsfonts}
\usepackage{amssymb}
\usepackage{graphicx}
\usepackage{cite}
\usepackage{algorithm2e}
\usepackage{array}
\usepackage[tight]{subfigure}
\usepackage{enumerate}

\IEEEoverridecommandlockouts
\begin{document}
\title{Mitigating Cascading Failures in Interdependent Power Grids and Communication Networks\thanks{This work was supported by DTRA grant HDTRA1-13-1-0021.}}

\author{\IEEEauthorblockN{Marzieh Parandehgheibi \IEEEauthorrefmark{1}, Eytan Modiano \IEEEauthorrefmark{1} and David Hay \IEEEauthorrefmark{2}}
\IEEEauthorblockA{\IEEEauthorrefmark{1}Laboratory for Information and Decision Systems, Massachusetts Institute of Technology, Cambridge, MA, USA\\
\IEEEauthorrefmark{2}School of Engineering and Computer Science, Hebrew University, Jerusalem 91904, Israel\\}}

\maketitle 

\begin{abstract}
In this paper, we study the interdependency between the power grid and the communication network used to control the grid. A communication node depends on the power grid in order to receive power for operation, and a power node depends on the communication network in order to receive control signals for safe operation. We demonstrate that these dependencies can lead to cascading failures, and it is essential to consider the power flow equations for studying the behavior of such interdependent networks. We propose a two-phase control policy to mitigate the cascade of failures. In the first phase, our control policy finds the non-avoidable failures that occur due to physical disconnection. In the second phase, our algorithm redistributes the power so that all the connected communication nodes have enough power for operation and no power lines overload. We perform a sensitivity analysis to evaluate the performance of our control policy, and show that our control policy achieves close to optimal yield for many scenarios. This analysis can help design robust interdependent grids and associated control policies.
\end{abstract}

\IEEEpeerreviewmaketitle

\section{Introduction}\label{Intro_sec}

One of the main challenges for sustainability of future power grids is the increased variability and uncertainty caused by integrating renewable sources into the grid. In order to address this challenge, the future grid must become equipped with real-time monitoring and be controlled with fast and efficient control algorithms \cite{IEC2012:Online}.

The monitoring and control of today's power grid relies on a Supervisory Control and Data Acquisition (SCADA) system. One of the main control operations is the Automatic Generation Control (AGC) which is used to match power supply with demand in the grid through frequency control.  This is done both at the local (generator) level, and the wide-area level. AGC systems rely on communications in order to disseminate control information, and the lack of communications, or even delay in communications can cause AGC systems to malfunction and fail, leading to wide-scale power outages \cite{Murty:Online,Green2014:Online,NERC2011:Online,GeneratorProtection:Online}. 

In August 2003, lack of real-time monitoring and rapid control decisions for mitigating failures led to a catastrophic blackout which affected 50 million people in Northeast America. According to the final report of the 2003 blackout \cite{FinalReport2004:Online}, this event started with the loss of transmission lines in Ohio due to inadequate tree trimming. However, the operators did not realize these failures due to insufficient monitoring; thus, no remedial action was taken at that time. In the subsequent hour, several transmission lines and generators tripped due to overheating of power lines and local protections\footnote{Local protections are systems that trip the generator when abnormal changes such as over/under frequency occur in the grid.} in generators.  These initial failures triggered a very fast cascade, which occurred in less than 5 minutes and led to a full blackout in the Northeast United States and parts of Canada. The reports in \cite{NERC2003:Online} and \cite{NY-ISO2004:Online} indicate that the reason for tripping of many generators and transmission lines was power imbalance in the control areas and lack of communication between the operators for mitigating the failures. It is thus essential to design a communication network together with control policies that facilitate widespread monitoring of the power grid, and enables the power grid to react to rapid changes and unexpected failures in the network. 

Moreover, for cost and sustainability considerations, the communication equipment often receives the power for operation directly from the power grid. However, this creates a strong interdependency between the two networks, where the operation of the power grid is dependent on receiving control signals from the communication network, and the operation of the communication network is dependent on receiving power from the power grid. 

The concept of interdependencies between infrastructures was first introduced in \cite{Rinaldi2001}. In \cite{Rosato2008}, Rosato \textit{et al.} studied the impact of failures in the power grid on the performance of communication networks. In 2010, Buldyrev et al. \cite{Buldyrev2010} presented a model for analyzing the robustness of interdependent random networks and investigated asymptotic connectivity to a ``giant component". They showed that interdependent networks are more vulnerable to failures than individual networks in isolation. The authors in \cite{Parandehgheibi2013} modified the model of \cite{Buldyrev2010} to account for connectivity to generators and control centers, and studied connectivity in the non-asymptotic regime.

\begin{figure}[ht]
\centering
\subfigure[Random Erdos-Renyi Graph - The size of largest component is close to the size of network for small sizes of failure, and there is a drop at a failure rate of about 50\%.]
{\label{fig:InetrdepNet}\includegraphics[scale=0.065]{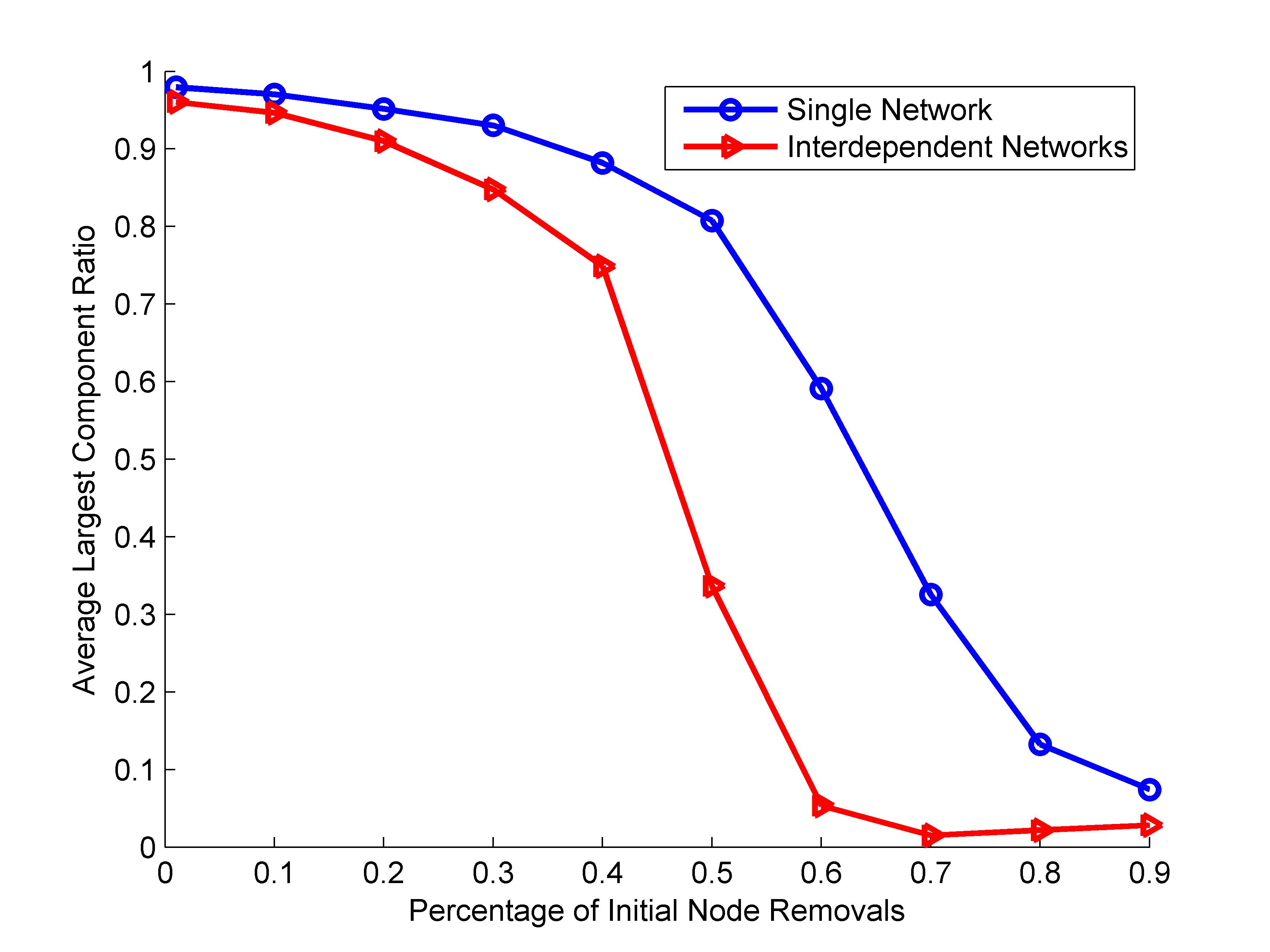}}
\subfigure[Random Power Grid - No large component exists for any size of failure.]
{\label{fig:InterdepPower}\includegraphics[scale=0.085]{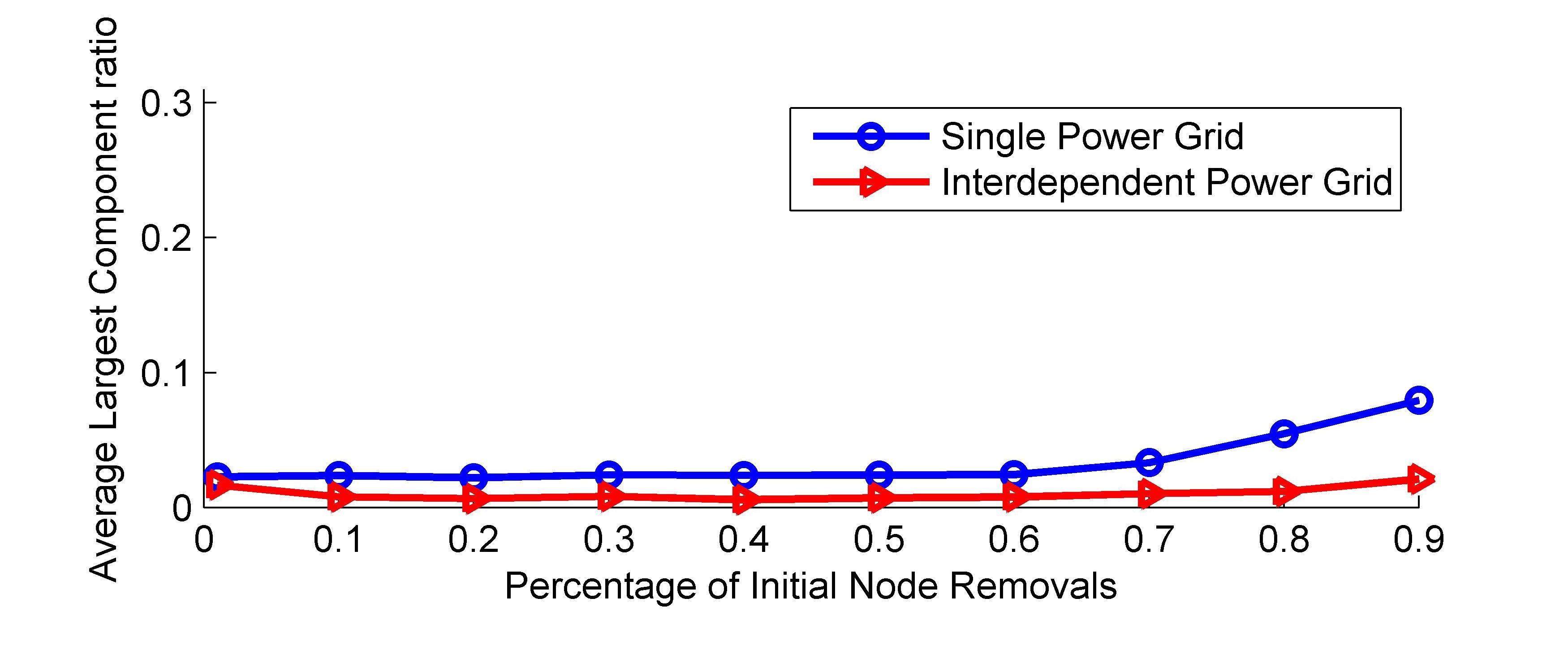}}                      
\label{Nature_comparison}
\caption{Ratio of largest component to the number of remaining nodes for different sizes of initial failures; Each network has 500 nodes with expected degree 4. We randomly selected $1/5^{th}$ of the nodes in the power grid and communication network as generators and control centers, respectively, and there is one-to-one interdependency between the two networks. For the power grid, we considered unit reactance for all power lines, and attributed a random amount of power in the range $[1000,2000]$ to all generators and loads.}
\vspace{-2mm}
\end{figure}

Figure \ref{fig:InetrdepNet} shows the impact of failures on interdependent networks by considering two random Erdos-Renyi graphs and a one-to-one interdependency between power and communication nodes. The result of our analysis shows that the interdependent networks are more vulnerable to failures than isolated networks, and there is a notable drop in the size of the largest connected component when the percentage of initial node removals is more than 50\% of total nodes. These results are in fact similar to those obtained by Buldyrev's model in \cite{Buldyrev2010}. 

However, in a power grid the flows are driven by Kirchoff's laws, and cannot be described by a network flow model. Thus, when a failure occurs in a power grid, the power flow is redistributed on the the rest of the network and some elements could overload and fail, leading to ``Cascading Failures". Since this behavior is not captured in the abstract models of \cite{Buldyrev2010,Parandehgheibi2013}, we generated random power grids and implemented the model of cascading failures from \cite{Bernstein2012} in conjunction with the model introduced in \cite{Parandehgheibi2013}. As can be seen from Figure \ref{fig:InterdepPower} when taking power failure cascades into account, there exists no large component for any size of initial failure. Thus, it is critical to consider the actual power flow in analyzing the behavior of the power grid.

\begin{figure}[ht]
	\centering
	\includegraphics[scale=0.085]{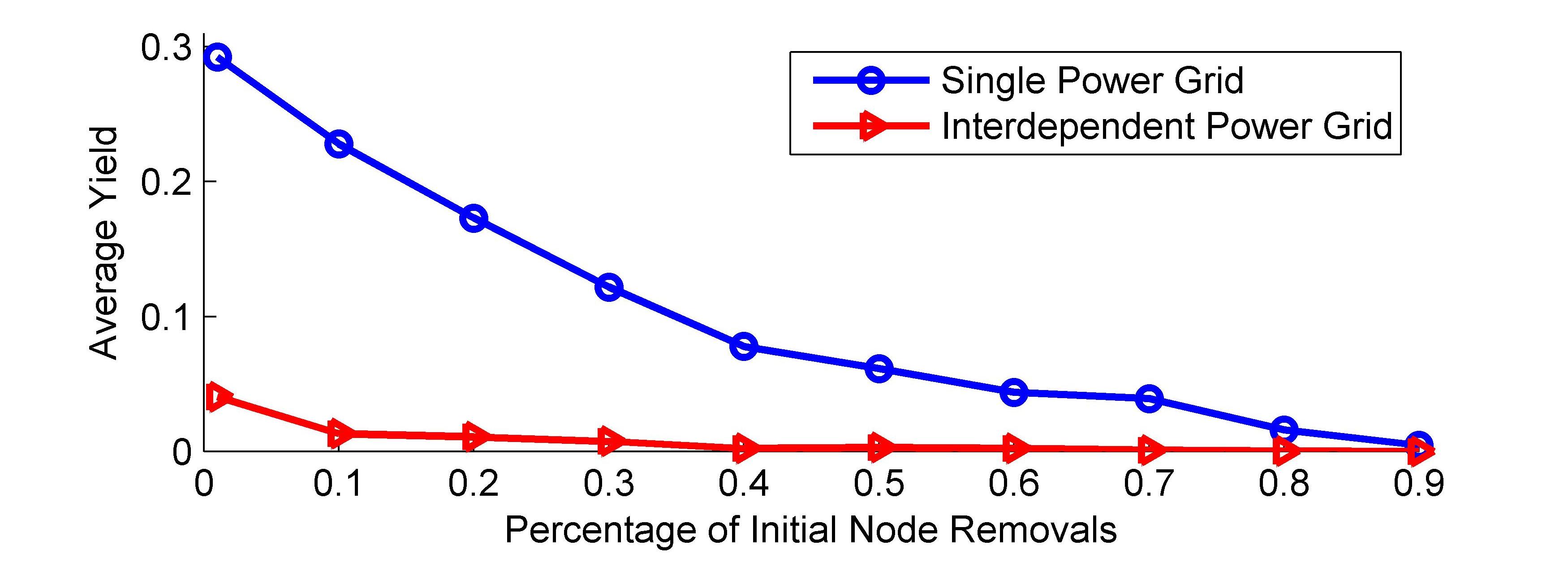}	
	\caption{Comparing Yield in single and interdependent Power Grid after cascading failures. Increasing the size of failure results in a smaller yield.}
	\label{Yield_Single_Interdep}
	\vspace{-2mm}
\end{figure}

Figure \ref{fig:InterdepPower} focused on connectivity. However, in power grids, the metric of interest is yield, which is the percentage of served load, not the the size of largest component. Figure \ref{Yield_Single_Interdep} shows that the yield in an interdependent power grid is much smaller than the yield in an isolated power grid.

As was discussed previously, it is necessary to design a communication network intertwined with the power grid in order to provide real-time monitoring and control for the grid. Therefore, a proper analysis of interdependent networks should account for the availability of control schemes that can mitigate cascading failures. In this paper, we propose a new load shedding scheme to control the cascade of failures both inside and between the networks. To the best of our knowledge, this paper is the first attempt to design control policies for mitigating failures in interdependent networks.

The rest of this paper is organized as follows. We explain the model of interdependent power grid and communication network in Section \ref{Model_Sec}. In Section \ref{Control_Sec}, we present a simple control policy followed by a load control policy that mitigate the cascading failures in interdependent networks in one stage. Finally, we perform sensitivity analysis for our load control policy in Section \ref{sensitivity_sec}and conclude in Section \ref{Conclusion_Sec}.

\section{Model}\label{Model_Sec}
An interdependent network consists of three subnetworks: Power grid, communication network and interdependency network. In the following, we will explain the model of each subnetwork.
\subsection{Power Grid}
The power grid can be modeled as a graph $G_P = (V_P, E_P)$ where $V_P$ and $E_P$ are power nodes and lines, respectively. There are three types of power nodes in a grid: Generators that generate power, Loads that consume power and Substations that neither generate nor consume power. 
The flow in power lines cannot be controlled manually; instead, it is determined based on the principles of electricity. In order to analyze the behavior of the power grid, we use the well-known DC power flow model, explained in equation \ref{DCFlow}, that has been widely used in the literature (see \cite{Bergen1999} for a survey on the power flow models).

Let $P$ be a $|V_P| \times 1$ vector such that $P_i$ denotes the power injection at power node $i \in V_P$. Let $A \in R^{|V_P| \times |E_P|}$ be the adjacency matrix where $A_{ij}=1$ if link $j$ starts from node $i$, $A_{ij}=-1$ if link $j$ ends in node $i$ and $A_{ij}=0$ otherwise. Moreover, let $X \in R^{|E_P| \times |E_P|}$ be the reactance matrix associated to the power grid where $X_{ii}$ denotes the reactance of $i^{th}$ power line and $X_{ij}=0$ for $j \neq i$. Let $f \in R^{|E_P| \times 1}$ be the vector of power flows in transmission lines and $\theta \in R^{|V_P| \times 1}$ denote the phases at all power nodes. A DC power flow can be modeled as follows.

\begin{subequations}
\vspace{-5mm}
\begin{align}
		 \quad & Af=P \label{PowFlowConst0} \\
		 \quad & A^T \theta=Xf \label{ReactConst0} 
\end{align}
\vspace{-5mm}
\label{DCFlow}
\end{subequations}

Constraint \ref{PowFlowConst0} is a network flow constraint which guarantees that power at every node is balanced. In addition, constraint \ref{ReactConst0} replicates Ohm's law where the amount of power flowing in a power line is equal to the difference in phase angles $\theta_i$ and $\theta_j$ divided by the reactance of line $(i,j)$.

When a power node or line fails, its load is shifted to other elements of the grid. During this process, the flow in one or more lines may be pushed beyond their capacity which leads to the failure of the overloaded lines. Similarly, failure of these lines redistributes power and may lead to further ``Cascading Failures". 

The cascade of failures in the power grid is a very complex phenomena, and several models have been introduced for explaining the behavior of cascading failures (see for example \cite{Bernstein2012,Chen2005,Dobson2003,Dobson2004}). In this paper, we will use the deterministic model explained in \cite{Bernstein2012}. In this model, each power line is associated with a capacity which is considered to be a factor of safety ($FoS$) typically set to $1.2$ times the amount of flow on that line. When a failure occurs, the power will be redistributed to the rest of the grid and the lines with flow more than their capacity will fail. The cascading model can be explained using the following steps. 

\begin{enumerate}
	\item Balance the power in the grid; i.e. if the grid is overloaded, decrease the amount of power at all loads uniformly to match the generation and if the grid is underloaded, decrease the amount of power at all generators to match the load.
	\item Resolve the DC power flow model in equations \ref{DCFlow}.
	\item Remove all the overloaded power lines; i.e. $f > f^{max}$. 
	\item If there is no overloaded lines in step 3, the cascade ends. Otherwise, repeat the four steps.
\end{enumerate}

\subsection{Communication Network}
The communication and control network can be modeled as a graph $G_C = (V_C, E_C)$ where $V_C$ and $E_C$ are communication nodes and links, respectively. There are two types of communication nodes: routers that are responsible for transmitting information, and control centers that are responsible for making control decisions.

In order to have a fully monitored and controlled power grid, every power node is equipped with a communication node (router). These nodes receive information from the power nodes and relay it to the control center through other routers. The control center makes the control decisions and sends them back to the routers located at power nodes. In our model, when a communication node fails, all the communication nodes that become disconnected from the control centers can no longer function. 

\subsection{Interdependency}
\subsubsection{Dependency of Communication Network on the Power Grid}
\begin{figure}[ht]
\centering
\includegraphics[scale=0.55]{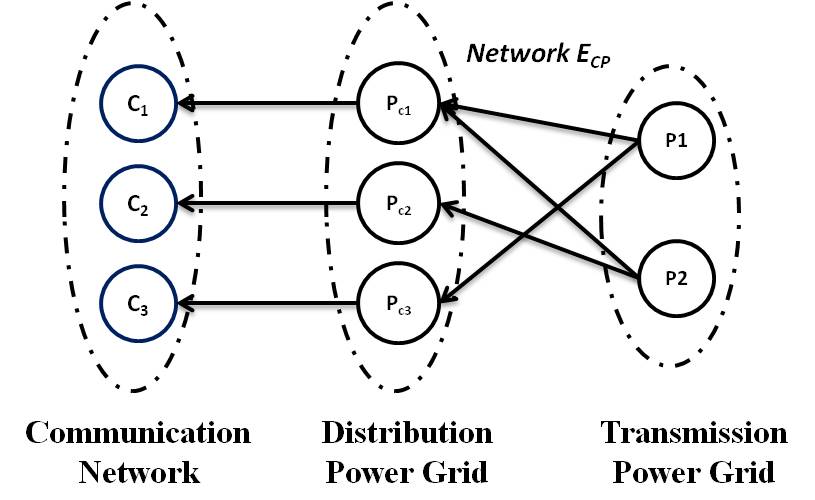}                   
\caption{Modeling Dependency of Communication Network on Power Grid}
\vspace{-3mm}
\label{fig:model}
\end{figure}

The communication nodes receive the power for their operation from power grid. In order to model this dependency, we associate each communication node $C_j$ with a load $P_{C_j}$ that is connected to the power grid (Figure \ref{fig:model}). Let $P_C^{req}$ be the required amount of power for operation of the communication node. Thus, communication node $C_j$ operates if $P_{C_j} \geq P_C^{req}$ and it fails otherwise.

In our model, the loads associated to the communication nodes are located in the distribution system and multiple communication nodes can receive power from one power node in the transmission system (Figure \ref{fig:model}). We assume that the communication and control nodes have the highest priority in the distribution system; thus, they will receive power as long as the power nodes have sufficient power to meet their demand. We model this part using network flow equations, where the sources are the loads $P_i$ in the power grid and the destinations are loads $P_{C_j}$ located at communication nodes. 

\subsubsection{Dependency of Power Grid on the Communication Network}
Next, we model the impact of loss of communication on the operation of the power grid. As explained in the introduction, AGCs control the operation of generators by setting the amount of power they should generate. If a generator becomes disconnected from the controller, the local controller tries to adjust the generation within a small range of changes in frequency. When the power grid is under stress (e.g. due to failures in the grid), power imbalance can lead to rapid frequency changes; in which case local protection schemes will be activated and trip the generators\cite{FinalReport2004:Online,NERC2011:Online,GeneratorProtection:Online}. Similarly, if a substation loses its control, then the relays cannot be accessed remotely and when the system is under stress, transmission lines will be overloaded and trip. Based on the 2003 blackout report, the large power imbalance in the system and lack of fast control and communication led to tripping of many transmission lines and generators \cite{FinalReport2004:Online}.

In this paper, we analyze the cascade of failures in the power grid when the system is under stress. We say that if a power node loses its correspondent communication and control node(s), it cannot be controlled and fails. This is a deterministic model that can be extended to a probabilistic model where the power node fails randomly with some probability. 

In the next section, we will propose control schemes that mitigate cascade of failures by shedding loads and re-dispatching generators. In our analysis, we do not study the transient behaviors of the grid after applying control decisions. Instead, we assume that due to a wide-area control implemented by the communication and control network, all the power nodes are aware of the transient changes in the system and local protections do not activate. This is essential as in the 2003 blackout, many generators tripped due to fluctuations resulting from intentional load sheddings \cite{FinalReport2004:Online}.

\section{Control Policies}\label{Control_Sec}
\vspace{1mm}
\subsection{Simple Load Shedding Mechanism}\label{Simple_Sec}
\vspace{1mm}
In this Section, we apply a simple load shedding control scheme in order to mitigate failures inside the power grid. This control scheme changes the power injection at power nodes so that the total power in the grid is balanced and the flow in transmission lines is below their capacity; thus, no failure cascades in the power grid. Different versions of this algorithm exist in the literature (see for example \cite{Koch2010, Bienstock2011}). The simple mitigation policy can be expressed in terms of the linear programming formulation in equation \ref{SimpleLoadShed_eqn}. Notice that notation ``updated" indicates that the power grid and communication network have been updated after initial failure. Let $P^{old}, P^{new} \in R^{|V^{updated}_P| \times 1}$ denote the power injections at power nodes before and after applying the simple mitigation policy. Moreover, let vector $f^{max}$ denote the capacity of power lines.

The objective function \ref{MinLoadShed1} is minimizing the total change in the power. Constraints \ref{maximum_gen1} and \ref{maximum_load1} enforce that the only possible controls are to shed loads and reduce power at generators. This is due to the fact that generators can ramp down much faster than they can ramp up. Since this control decision should be applied very rapidly in order to keep the network stable, we only allow ramping down; i.e. decreasing generation. Moreover, we assume there is no minimum threshold on the amount of power generation or consumption.

\begin{subequations}
\vspace{-3mm}
\begin{align}
\mbox{minimize} \quad & e^T(|P^{new} - P^{old}|) \label{MinLoadShed1} \\
\mbox{subject to} \quad & A^{updated}f=P^{new} \label{PowFlowConst1} \\
									\quad & (A^{updated})^T \theta=Xf \label{ReactConst1} \\
									\quad & f \leq f^{max}  \label{CapacityConst1} \\
									\quad & 0 \leq P_{i}^{new} \leq P_{i}^{old} \quad \forall i\in V_{P,gen}^{updated} \label{maximum_gen1}	\\
									\quad & P_{i}^{old} \leq P_{i}^{new} \leq 0 \quad \forall i\in V_{P,load}^{updated} \label{maximum_load1}						
\end{align}
\label{SimpleLoadShed_eqn}
\vspace{-3mm}
\end{subequations}

We apply this mitigation policy to the random interdependent power grid generated in Section \ref{Intro_sec}. The only difference is that communication nodes receive power only from loads; thus, it is not a fully one-to-one interdependent topology. However, we try to create as many one-to-one interdependencies as possible; i.e a load is dependent on the communication node that it provides power for. Previous studies have shown that one-to-one interdependent networks are more robust to failures \cite{Parandehgheibi2013}. We observed that although applying this control policy can mitigate failures inside power grid, the failures still cascade between the communication network and the power grid. Thus, we apply the control algorithm iteratively until no further failures occur. Clearly, the yield in any interdependent topology would be upperbounded by the yield in an isolated power grid. We use this upperbound to examine the performance of our control scheme.

\begin{figure}[ht]
	\centering
	\includegraphics[scale=0.054]{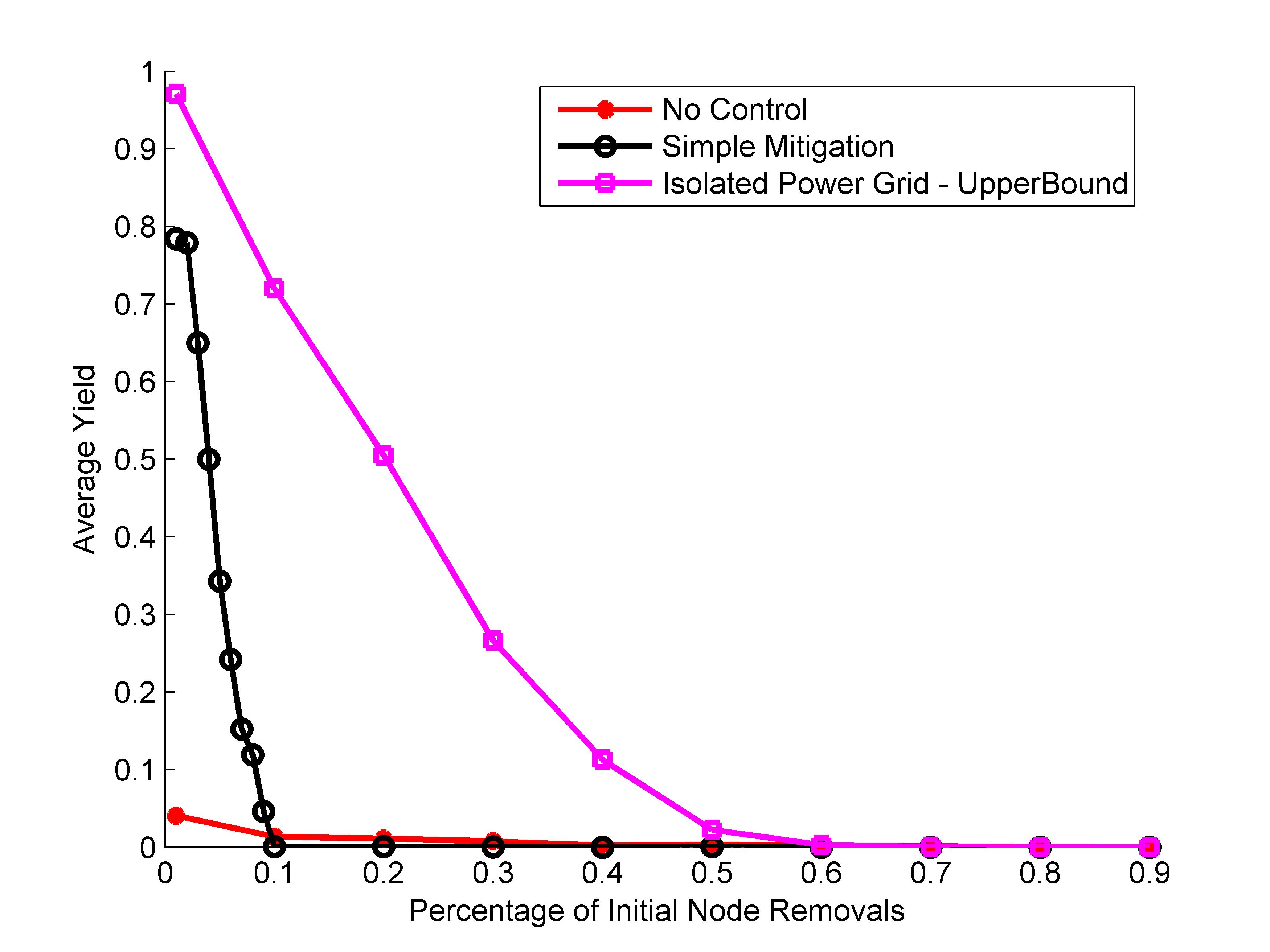}
	\caption{Applying Simple Mitigation Policy to Interdependent Power Grids for controlling Cascading failures. Here, the total power required by communication network is $10^{-4}$ times total load in the power grid.}
	\label{fig:SimpleLoadShedding}
	\vspace{-2mm}
\end{figure}

Figure \ref{fig:SimpleLoadShedding} shows the yield after applying the simple mitigation policy. It can be seen that although the control policy has improved the yield (when the failure rate is small), there is a dramatic drop at the beginning and the control policy can not survive any failure larger than 10\% of the network. This is due to the fact that loss of loads leads to the loss of communication nodes that are used to control the generators and thus, the generators fail. Therefore, it is much harder to mitigate the cascading failures in interdependent networks. This simple policy is meant to demonstrate that even simple controls can reduce cascades. Inspired by this observation, we develop a control scheme that aims to keep communication nodes operating.

\subsection{Load Control Mitigation Policy}\label{General_Sec}

It was seen in Section \ref{Simple_Sec} that a simple mitigation policy cannot mitigate failures in interdependent networks, as failures in the power grid propagate to the communication network and cause additional failures both inside the communication network and in the power grid. In order to avoid such propagations, we propose a novel control policy that consists of two phases: in the first phase, it predicts the non-avoidable failures in the power grid and the communication network and removes these nodes from the network. In the second phase, it changes the power injection at power nodes so that (1) power in all transmission lines is below their capacity and (2) all the remaining communication nodes keep operating; i.e. $P_{C_j} \geq P_C^{req}$. This guarantees that no further failures occur in the power grid and that the failures do not propagate to the communication network. Thus, the cascade of failures will be mitigated in one stage. In the following, we explain these phases in more details.

\subsubsection{Phase I}
In this phase, we ignore the power flows in the power grid and find the nodes that their failure cannot be avoided by changing the power injection at nodes due to loss of connectivity. Algorithm I describes how to find such failures in polynomial time. %Let $N$ be the total number of nodes in the network; then, the cascade algorithm terminates after at most $N$ iterations, because the algorithm should remove at least one node at each iteration. Moreover, in each iteration, at most $N$ nodes are checked; thus, the total run time is $O(N^2)$.

\setlength{\algomargin}{-0.3em} 
\begin{algorithm}

 \SetKwData{Left}{left}
  \SetKwData{Up}{up}
  \SetKwFunction{FindCompress}{FindCompress}
  \SetKwInOut{Input}{input}
  \SetKwInOut{Output}{output}
 
  \Input{Topology of interdependent network and the set of initial node removals}
\Repeat{No node can be removed}{
		\begin{enumerate}
			\item{For every power node $i$, check if there exists a path \\from a generator to node $i$ and it receives an incoming edge from the communication network;}
			\item{For every communication node $j$, check if there exists\\ a path from a control center to node $j$ and it receives \\an incoming edge from power grid;}
			\item{Remove all the nodes that do not satisfy the properties in steps 1 and 2;}
			\item{Remove all isolated generators;}
			\item{Remove all the links connected to the removed nodes.}
		\end{enumerate}
				}
	\Output{Set of all removed nodes \vspace{2 mm}}
\caption{Cascade Algorithm}
\end{algorithm}

\subsubsection{Phase II}
In this phase, our objective is to find a set of feasible power injections so that the minimum amount of load is shed and no control node fails due to loss of power. Let $E_{CP}$ denote the adjacency matrix modeling the dependency of communication network on loads. Let vector $h$ denote the amount of power flowing from loads located in the power transmission grid $P_i$ to loads located in the power distribution grid $P_{C_j}$ that support the communication network (See Figure \ref{fig:model}). Moreover, let $b$ be a two part power vector. The first part represents the amount of power injection at loads in the transmission grid with positive sign as these are source nodes; i.e. $-P_i$ since loads $P_i$ are originally modeled with negative values. Similarly, the second part represents the amount of power injection at loads in the distribution grid with negative sign as these are destination nodes; i.e. $P_{C_j}$ since loads $P_{C_j}$ are originally modeled with negative values. Moreover, notice that notation ``updated" indicates that the power grid and communication network have been updated by removing the nodes that fail in Phase I.

\begin{subequations}
\vspace{-4mm}
\begin{align}
\mbox{minimize} \quad & e^T(|P^{new} - P^{old}|) \label{MinLoadShed2} \\
\mbox{subject to} \quad & A^{updated}f=P^{new} \label{PowFlowConst2} \\
									\quad & (A^{updated})^T \theta=Xf \label{ReactConst2} \\
									\quad & f \leq f^{max}, \quad \forall (i,j)\in E^{updated}_P  \label{CapacityConst2} \\
									\quad & 0 \leq P_{i}^{new} \leq P_{i}^{old} \quad \forall i\in V_{P,gen}^{updated} \label{maximum_gen2}	\\
									\quad & P_{i}^{old} \leq P_{i}^{new} \leq 0 \quad \forall i\in V_{P,load}^{updated} \label{maximum_load2}	\\				 	
									\quad & P_{C_j} \leq -P_C^{req} \quad \forall j \in V^{updated}_C  \label{Guaranteelivecontrol}\\
									\quad & \sum E^{updated}_{CP}h=b \label{DistConst}\\
									\quad & h \geq 0  \label{DistFlowConst}
\end{align}
\vspace{-4mm}
\label{ControlPolicyEqn}
\end{subequations}

Constraint \ref{Guaranteelivecontrol} guarantees that every remaining communication node receives the minimum amount of power required for operation. Constraint \ref{DistConst} models the power flowing from the power grid to the communication network with a network flow model. Constraint \ref{DistFlowConst} shows that the direction of power flow is from power nodes to communication nodes. The combination of these three constraints changes the power injection at power nodes so that the communication nodes remaining from Phase I will continue operating. The rest of constraints are similar to the Simple mitigation policy and set the power injections so that no transmission line is overloaded. Similarly, the objective function is minimizing the total amount of load shedding. Note that this ILP may be infeasible. In such cases the yield would be zero, showing that our control policy is not capable of controlling the failures in the network.

\begin{figure}[ht]
\centering
\includegraphics[scale=0.054]{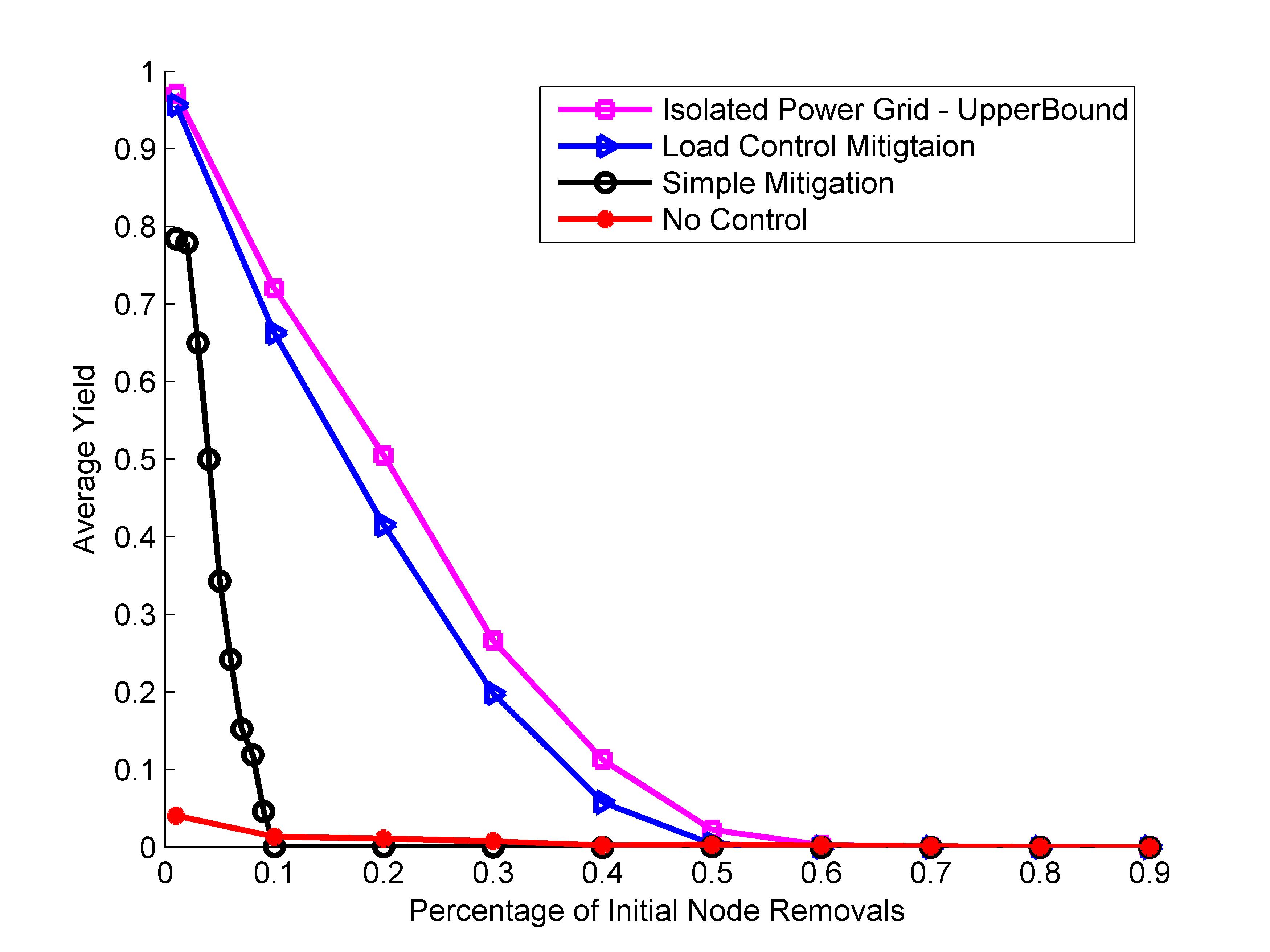}                 
\caption{Comparing performances of control policies}
\label{fig:EpsPerformance}   
\vspace{-2mm}
\end{figure}

Figure \ref{fig:EpsPerformance} compares the performance of load control mitigation and simple mitigation policies. It can be seen that the yield after applying the load control policy is improved with respect to the simple mitigation policy and it is very close to the upperbound. 

\section{Sensitivity Analysis}\label{sensitivity_sec}
We analyze the performance of our control policy with respect to changes in the communication network and the interdependency between the power grid and communication network. The parameters we study are the amount of power that communication nodes require ($P_{C_j}$), size of communication network, the average number of power nodes supporting each communication node, namely ``Communication Interdependence Degree" and finally, the average number of communication nodes supporting each power node, namely ``Power Interdependence Degree". 

We generate 30 random networks and test their feasibility by applying our control policy. If the entire network fails in the first phase, the network is not feasible; i.e. no control policy can survive it. We average the yield found by our control policy over the feasible networks.

\begin{figure}[ht]
\centering
\includegraphics[scale=0.054]{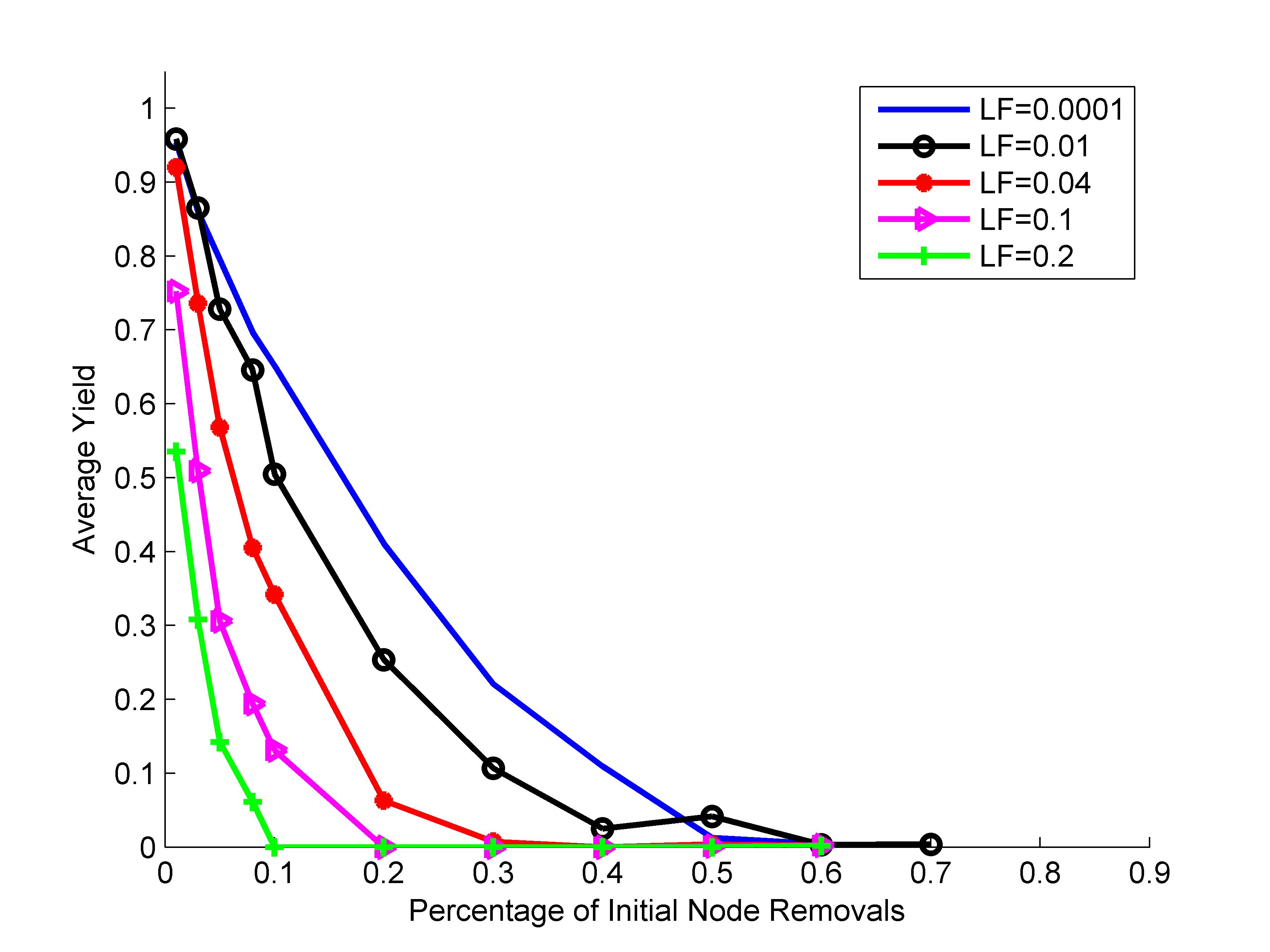}
\caption{Impact of Load Factor on the Yield; 500 communication nodes; 20\% of nodes are randomly selected as control centers}
\label{fig:yield_LF}
\vspace{-2mm}
\end{figure}

We define the ``Load Factor" as the ratio of power required by the communication network to the total load in the power grid. In the previous simulations LF was set to be $10^{-4}$. Figure \ref{fig:yield_LF} shows that by increasing LF, the yield decreases as it is harder to provide a larger amount of power for all loads supporting communication network. 

\begin{figure}[ht]
\centering
\includegraphics[scale=0.054]{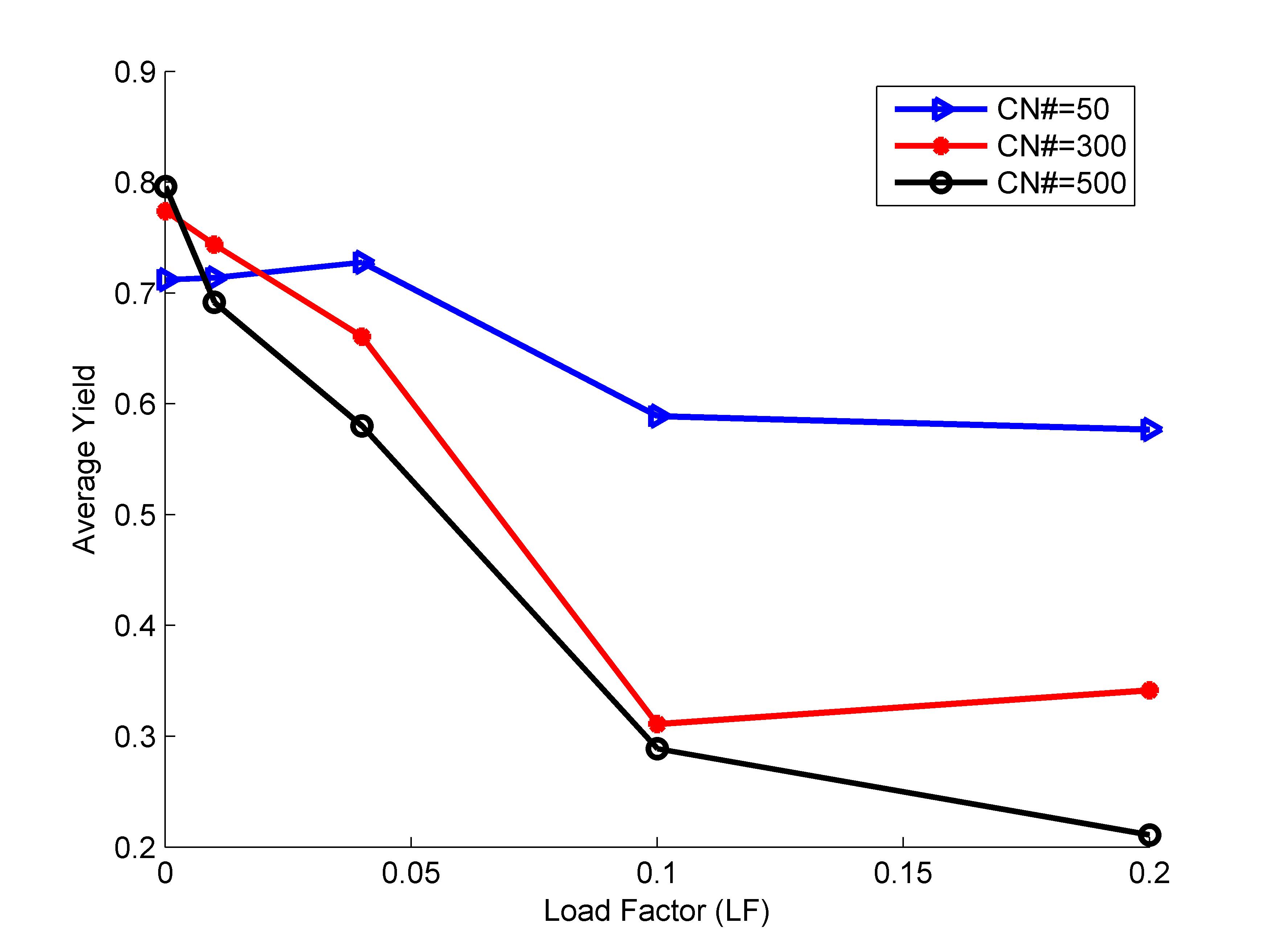}
\caption{Impact of Number of Communication Nodes on the Yield; 20\% of communication nodes are randomly selected as control centers; initial removal=5\%}
\label{fig:yield_ComNode}
\vspace{-2mm}
\end{figure}

We analyze the performance of our policy with respect to the size of communication network; i.e. number of communication nodes. Figure \ref{fig:yield_ComNode} shows that for small values of LF, the larger networks have higher yield. However, by increasing LF, the yield of larger networks decreases more; thus, the smaller networks perform better for large LF. 

\begin{figure}[ht]
\centering
\includegraphics[scale=0.054]{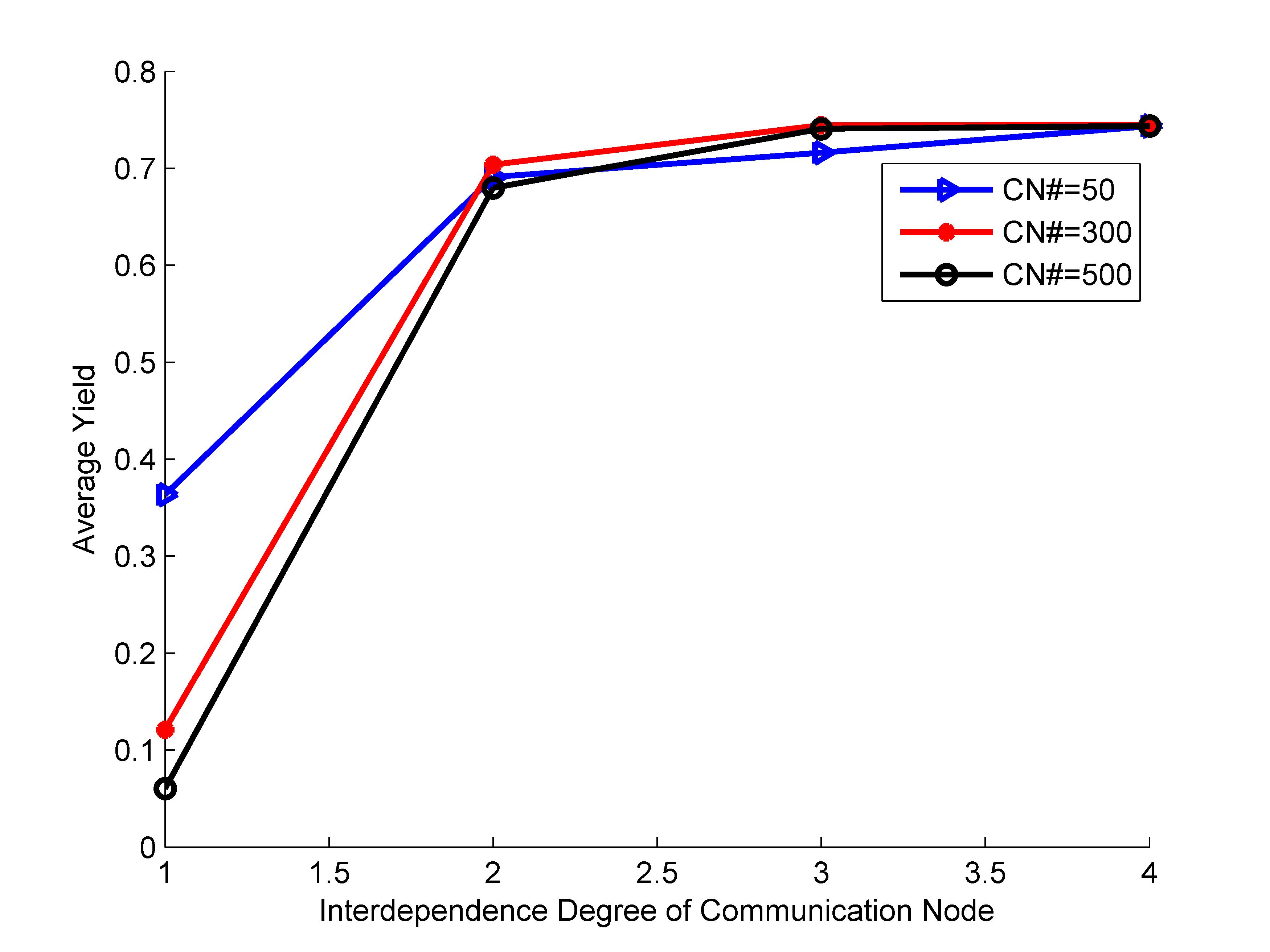}
\caption{Impact of Communication Interdependence Degree on the Yield; initial removal=10\%; load factor=0.1}
\label{fig:yield_ComDegree}
\vspace{-5mm}
\end{figure}

The next parameter that we study is the average number of power nodes that support every communication network (interdependence degree). Figure \ref{fig:yield_ComDegree} shows that the average yield increases by increasing the interdependence degree. Moreover, it shows that for large enough degree (here degree of 4) networks with different sizes have similar yield.

\begin{figure}[ht]
\centering
\includegraphics[scale=0.054]{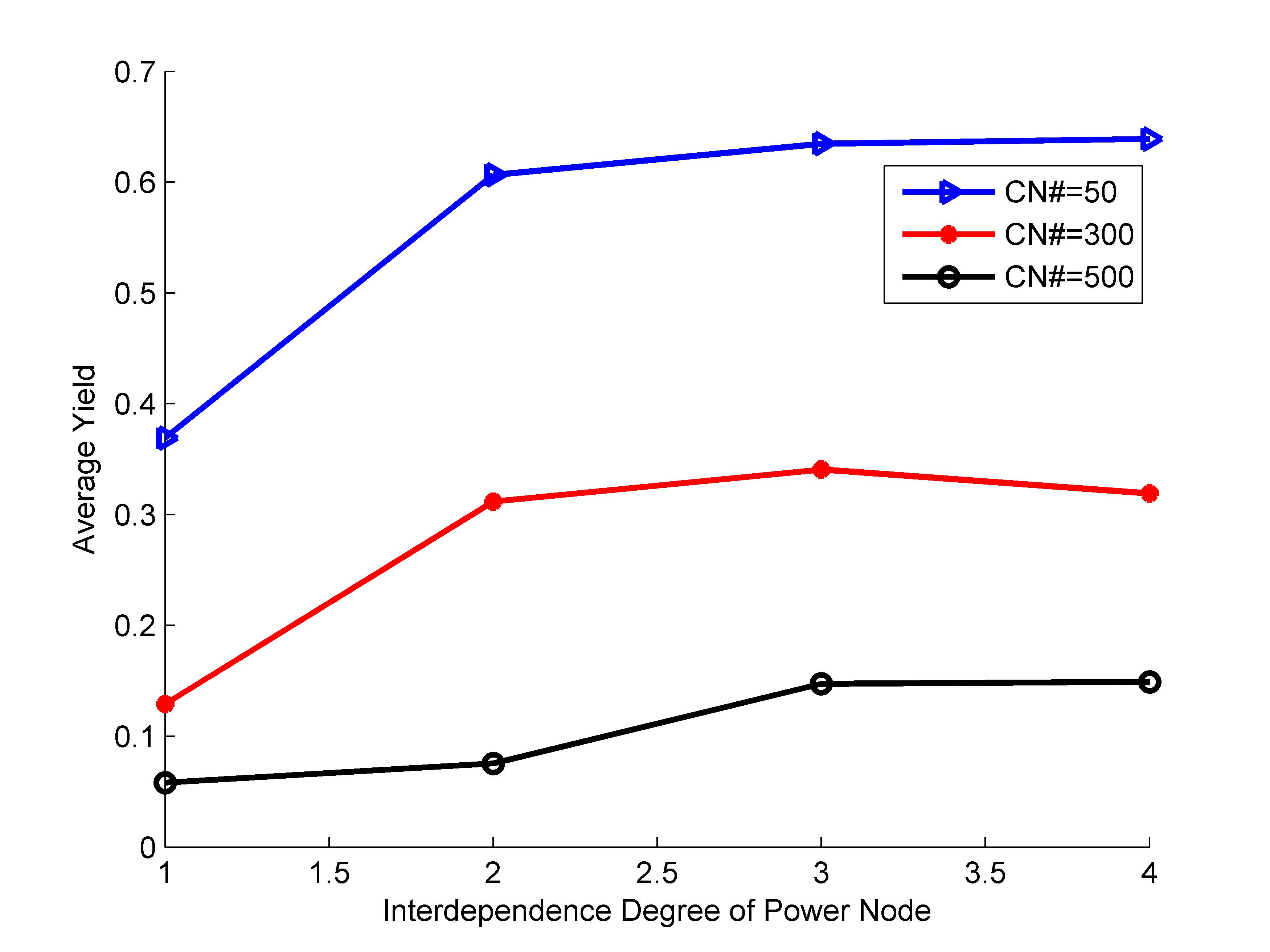}
\caption{Impact of Power Interdependence Degree on the Yield; initial removal=10\%; load factor=0.1}
\vspace{-5mm}
\label{fig:yield_PowDegree}
\end{figure}

Finally, we investigate the impact of the average number of communication nodes supporting each power node. It can be seen from Figure \ref{fig:yield_PowDegree} that increasing degree has positive impact on the yield and feasibility; however, it is not as strong as the impact of communication interdependence degree. The reason is due to the structure of our control policy that tends to survive all of the communication network. Thus, if the control policy is feasible, the communication network remains operating, which results in the operation of the power nodes supported by these communication nodes. Therefore, in these scenarios, increasing the support for power grid cannot help. The small improvement that we see here is related to the reduction of failures due to disconnection. 

\section{Conclusion}\label{Conclusion_Sec}
In this paper, we showed that it is essential to consider the power flow equations for analyzing the behavior of interdependent power grid and communication networks. We argued that in order to analyze the robustness of interdependent networks, one should consider the control schemes for controlling cascading failures both inside and between the power grid and communication network. We proposed a new control scheme that mitigates failures in one stage and keeps the yield close to the maximum possible value. Our policy only allowed the failures due to disconnection from generators and control centers. Thus, a connectivity model can be used to describe the process of cascading failures in interdependent topologies. In addition, we tested the performance of our load control policy with respect to changes in several parameters such as load fraction (power needed by communication) and interdependence degree.

\bibliographystyle{IEEEtran}
\bibliography{biblio}

\end{document}